\documentclass[11pt]{amsart}
\usepackage{epsfig}
\usepackage{graphics}

\newtheorem{theorem}{Theorem}

\newtheorem{corollary}[theorem]{Corollary}

\theoremstyle{definition}

\newtheorem{conj}{Conjecture}

\theoremstyle{remark}

\def\mod{{\rm Mod}}

\begin{document}

\newenvironment{prooff}{\medskip \par \noindent {\it Proof}\ }{\hfill
$\square$ \medskip \par}
    \def\sqr#1#2{{\vcenter{\hrule height.#2pt
        \hbox{\vrule width.#2pt height#1pt \kern#1pt
            \vrule width.#2pt}\hrule height.#2pt}}}
    \def\square{\mathchoice\sqr67\sqr67\sqr{2.1}6\sqr{1.5}6}
\def\pf#1{\medskip \par \noindent {\it #1.}\ }
\def\endpf{\hfill $\square$ \medskip \par}
\def\demo#1{\medskip \par \noindent {\it #1.}\ }
\def\enddemo{\medskip \par}
\def\qed{~\hfill$\square$}

 \title[On a question of Brendle and Farb]
{On a question of Brendle and Farb }

\author[Mustafa Korkmaz]{Mustafa Korkmaz$\,^*$ }

 \address{Department of Mathematics, Middle East Technical University,
 06531 Ankara, Turkey} \email{korkmaz@arf.math.metu.edu.tr}

 \date{\today}
\thanks{* Supported by T\"UBA/GEB\.{I}P}
\begin{abstract}
In this note we prove that there is no constant $C$, depending on 
the genus of the surface, such that 
every element in the mapping class group can be written as a 
product of at most $C$ torsion elements, answering a question of 
T. E. Brendle and B. Farb in the negative.\\


\end{abstract}

 \maketitle

  \setcounter{secnumdepth}{1}
 \setcounter{section}{0}

\section{Introduction} 

In a recent paper~\cite{bf}, Tara E. Brendle and Benson Farb showed that the mapping class group
of a closed orientable surface of genus $g\geq 3$ is generated by three torsion elements,
and also by seven involutions. They asked whether there is a number
$C=C(g)$ such that every element in the mapping class group  can be written as a 
product of at most $C$ torsion elements. The purpose of this note is to give a negative
answer to this question. We deduce our result from the classification of torsion 
elements in the mapping class group and the fact that the mapping class group is 
not uniformly perfect, although it is perfect.

\section{Preliminaries}
Let $S$ be a closed orientable surface of genus $g$ and let $\mod_g$ denote
the mapping 
class group of $S$, the group of the isotopy classes of orientation-preserving 
diffeomorphisms $S\to S$. 

The following theorem is due to Kazuo Yokoyama~\cite{y}. In fact, Yokoyama
gives the number of conjugacy classes explicitly.
 
\begin{theorem}   \label{thm:periodic}
The number of conjugacy classes of torsion elements in the
mapping class group $\mod_g$ is finite.
 \end{theorem}

\begin{corollary}  \label{cor:periodic}
Suppose that $g\geq 3$.
There is a constant $T(g)$ such that every torsion element in the
mapping class group $\mod_g$ can be written as a product of at most $T(g)$
commutators.
 \end{corollary}
\begin{proof}
Choose a representative from each conjugacy classes and write them as a product of 
commutators. Let $T(g)$ be the maximum number of such commutators.
Since the conjugation of a commutator is again a commutator, the corollary
follows from Theorem~\ref{thm:periodic}.
\end{proof}

In the above corollary and its proof, we may replace ``commutator(s)" by ``Dehn twist(s)".

It is well known that the mapping class group $\mod_g$ is perfect for
$g\geq 3$, but it is not uniformly perfect.

\begin{theorem} $($\cite{k},\cite{ek}$)$ \label{thm:stable}
Suppose that $g\geq 3$.
There is no constant $C(g)$ such that every element in the
mapping class group $\mod_g$ can be written as a product of at most $C(g)$
commutators.
 \end{theorem}

\section{The results}

\begin{theorem}  \label{result}
Suppose that $g\geq 3$.
There is no constant $C(g)$ such that every element in the mapping class group
$\mod_g$ of a closed orientable surface of genus $g$ can be written as a product of at most $C(g)$ 
torsion elements.
 \end{theorem}
\begin{proof}
Suppose that there is such a constant $C(g)$. We then conclude from 
Corollary~\ref{cor:periodic} that every element in $\mod_g$ is 
a product of at most $C(g)T(g)$ commutators, 
which contradicts to Theorem~\ref{thm:stable}.
\end{proof}

Corresponding question can also be asked for Dehn twists; is there a 
constant $C(g)$ such that every element in $\mod_g$ is a product of 
at most $C(g)$ Dehn twists? This question, which was asked to the author
by Andr\'as Stipsicz, was the motivation for the paper~\cite{k}.
Since every Dehn twist can be written as a product of two commutators~\cite{ko},
the following theorem, which is implicit in \cite{k} and~\cite{ek},
follows easily from Theorem~\ref{thm:stable}. We note that it also holds 
for the mapping class group $\mod_2$.

\begin{theorem}  
Suppose that $g\geq 3$.
There is no constant $C(g)$ such that every element in the mapping class group
of a closed orientable surface can be written as a product of at most $C(g)$ 
Dehn twists.
 \end{theorem}

\section{Final comments}
Any element $f\in\mod_g$ can be written as a product of torsion
elements~\cite{bf}. Let us define $\tau_g(f)$ to be the minimum number of such torsion 
elements and call it the {\it torsion length} of $f$. 
Clearly, the sequence $\tau_g(f^n)$ is subadditive, 
that is, $\tau_g(f^{n+m})\leq \tau_g(f^n)+\tau_g(f^m)$ for all positive integers $n,m$.
Therefore, the limit $$\lim_{n\to \infty}\frac{\tau_g(f^n)}{n}$$ exists. Let us denote this limit by
$||f||_{\tau_g}$ and call it the {\it stable torsion length} of $f$. 

If $t_a$ denotes the Dehn twist about some nonseparating
simple closed curve $a$, then $t_a$ is a product of two torsion elements 
(see~\cite{bf}, Lemma $3$). Hence, $\tau_g(t_a)=2$. What is $\tau_g(t_a)$ if $a$
is separating? In general, what is $\tau_g(t_a^n)$? The fact
$\tau_g(t_a)=2$ gives the upper bound $||t_a||_{\tau_g}\leq 2$.

On the other hand, for any nontrivial simple closed curve $a$ 
the growth rate of $t_a$ is linear~\cite{k,flm}
and its stable commutator length is positive~\cite{k,ek} (see the references for 
definitions).  Inspired by these facts, the following conjecture seems reasonable.

\begin{conj}
The stable torsion length $||t_a||_{\tau_g}$ is positive.
\end{conj}

Clearly, if $f$ is a torsion element, then the torsion length of any power of $f$ is 
either one or zero. Therefore, its stable torsion length $||f||_{\tau_g}$ is zero. We end
by stating a stronger version of Conjecture~$1$.

\begin{conj}
If $||f||_{\tau_g}=0$, then $f$ is torsion. 
\end{conj}

\end{document}